\documentclass{aims}

\usepackage{txfonts}
\def\typeofarticle{Research article}
\def\currentvolume{5}
\def\currentissue{x}
\def\currentyear{2021}

\def\ppages{xxx--xxx}
\def\DOI{}
\def\Received{March 2021}
\def\Accepted{xxxx xxxx}
\def\Published{xxxx xxxx}

\usepackage{IEEEtrantools}

\usepackage{graphicx,amsmath,amsfonts,amsthm}

\numberwithin{equation}{section}

\newcommand{\norm}[1]{\left\lVert#1\right\rVert}
\newcommand{\nor}[1]{\parallel #1 \parallel}

\def\R{\mathbb{R}}

\def\bu{\textbf{u}}

\def\bn{\textbf{n}}

\def\bS{\textbf{S}}
\def\bs{\textbf{s}}

\def\Bila{\mathbf{B}}

\def\Ne{\mathcal{N}}

\def\P{\mathcal{P}}

\def\A{\mathcal{A}}

\def\T{\mathcal{T}}

\newcommand{\eequ}{\end{equation}}

\newcommand{\bequ}{\begin{equation}}

\newcommand{\re}[1]{\textcolor{black}{#1}}
\newcommand{\dt}[1]{\textcolor{black}{#1}}

\newcommand{\dtr}[1]{\textcolor{black}{#1}}

\begin{document}

\title{Nonlocal multiscale \re{ modelling  of} tumour-oncolytic viruses interactions within a \re{heterogeneous fibrous/non-fibrous} extracellular matrix}

\author{%
  Abdulhamed Alsisi\affil{1},
  and
  Raluca Eftimie\affil{2},
  and
  Dumitru Trucu\affil{1,}\corrauth
}

\shortauthors{the Author(s)}

\address{%
  \addr{\affilnum{1}}{Division of Mathematics, University of Dundee, Dundee DD1 4HN, United Kingdom}
  \addr{\affilnum{2}}{Laboratoire Mathematiques de Besan\c{c}on, UMR-CNRS 6623, Universit\'{e} de Bourgogne Franche-Comt\'{e}, 16 Route de Gray, Besan\c{c}on, France}
  }

\corraddr{trucu@maths.dundee.ac.uk}

\begin{abstract}
\re{In this study we investigate computationally tumour-oncolytic virus (OV) interactions that take place within a heterogeneous ExtraCellular Matrix (ECM).} The ECM is viewed as a mixture of two constitutive phases, namely a fibre phase and a non-fibre phase. \re{The multiscale mathematical model presented here focuses on the nonlocal cell-cell and cell-ECM interactions, and how these interactions might be impacted by the infection of cancer cells with the OV. At macroscale we track the kinetics of cancer cells, virus particles and the ECM. At microscale we track (i) the degradation of ECM by matrix degrading enzymes (MDEs) produced by cancer cells, which further influences the movement of tumour boundary; (ii) the re-arrangement of the microfibres that influences the re-arrangement of macrofibres (i.e., fibres at macroscale).}
\re{With the help of this new multiscale model, we investigate two questions: (i) whether the infected cancer cell fluxes are the result of local or non-local advection in response to ECM density; and (ii) what is the effect of ECM fibres on the the spatial spread of oncolytic viruses and the outcome of oncolytic virotherapy.}


\end{abstract}


\maketitle
\section{Introduction}
\label{Sect:Intro}

Cancer invasion is a complex \re{multiscale} phenomenon \re{that relies on the structure and composition of the extracellular matrix (ECM)}. The ECM is an important biological structure that serves as a platform for cellular communication, as well as providing support to surrounding cells and tissues, transducing mechanical signals, and functioning as adhesive substrate~\cite{Rozario2010}. \re{This matrix is a highly dynamic structure that controls most fundamental behaviours of cells: from proliferation, to adhesion, migration, differentiation and apoptosis~\cite{Yue2005_BiologyECM}. These processes are controlled through the interactions of cells with the components of the ECM: collagen, proteoglycans, elastin and cell-binding glycoproteins~\cite{Yue2005_BiologyECM}. The matrix components are continuously deposited, degraded or modified, and thus the ECM is continuously undergoing remodelling. This remodelling process (which involves matrix degrading enzymes that can be secreted by the cancer cells and normal cells) impacts also the evolution of cancer invasion. A fundamental process in the invasive potential of cancer cells is the cell-ECM adhesion (through adhesion molecules such as integrins)~\cite{Gkretsi2018_CellAdhesionECM}. In addition, the collective movement of cancer cells is the result of cell-cell adhesion (through adhesion molecules such as E-cadherins)~\cite{Gkretsi2018_CellAdhesionECM}. Understanding these cell-cell and cell-matrix adhesion processes is important not only for our understanding of the evolution of cancer and its invasion of the surrounding tissue, but also for our understanding of the efficacy of various anti-cancer therapies.}

An emerging effective anti-cancer treatment is the oncolytic virotherapy (OV-therapy). \re{The effectiveness of this treatment} lies on viruses selectively \re{infecting and destroying} malignant cancer cells without harming the surrounding healthy cells; see \cite{Fountzilas2017_OVreview,Kaufman2015,Pol2016,Russell2012}. \re{The structure and composition of the ECM influences the effectiveness of the OV-therapies, since the ECM components can form a physical barrier that traps the viral particles~\cite{Wojton2010_TumEnvir-OV}. In particular, experimental studies have shown that the collagen fibres play a very important role in inhibiting viral spread~\cite{Wojton2010_TumEnvir-OV}. }

\re{Due to the complexity of the tumour microenvironment, which makes it difficult to understand the interactions between the different components of this environment, mathematical models have been used over the last few decades to answer various questions about these interactions. The great majority of these models are single-scale models, which focus on spatial tumour invasion~\cite{Armstrong2006,Gerisch2008,Domschke2014a}, on tumour oncolytic therapies~\cite{CrivelliKimWares2012,Eftimie2011,Eftimie2016,GevertzWares2018_MathOV,Heidbuechel2020_MathModelOVReview,NowakM.A.MartinA.2000Vd:m,Wodarz2016}, or both ~\cite{Berg2019,Jacobsen2015,Malinzi2015,TaoWinkler2020_PDEMathModelOV,Wodarz2012}. More recently, various multi-scale mathematical models have been derived to reproduce and investigate biological processes that take place at different spatial scales~\cite{Alsisi2020,aalsisi21,Alzahrani2019,Alzahrani2020a,Paiva2009,Paiva2013,Trucu2013,Shuttleworth2019}. For example, \cite{Trucu2013} introduced a multi-scale moving boundary model for cancer invasion, which focused on the local interactions between cancer cells and the ECM, via matrix degrading enzymes (MDEs) that act at the micro-scale level of the invading tumour boundary. This study was further extended in~\cite{Shuttleworth2019} by considering also the non-local cell-cell and cell-ECM adhesive interactions. This study also considered a heterogeneous ECM population formed of fibres and non-fibrous sub-populations, and investigated the role of fibres (and their re-arrangement at macro-scale and micro-scale) in the evolution of a tumour population. Other local and nonlocal multiscale models were used in~\cite{Alzahrani2019,Alzahrani2020a,Alsisi2020,aalsisi21} to investigate the  interactions between oncolytic viruses and cancer cells. 
}

\begin{figure}[h]
    \centering
    \includegraphics[width=3.5in]{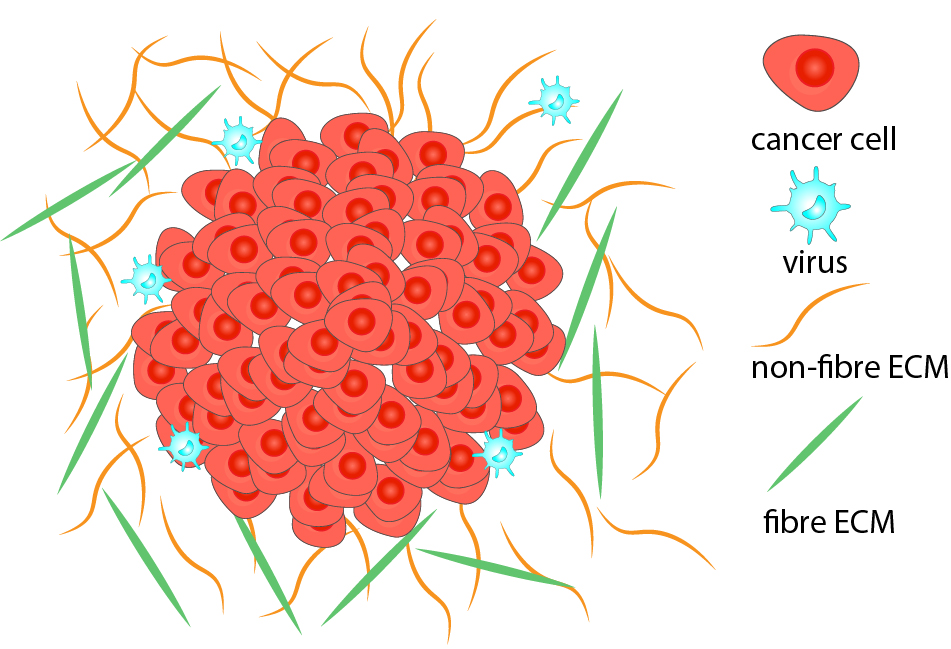}
    \caption{Caricature description of the interactions between oncolytic viruses (OVs), cancer cells, non-fibre ECM and fibre ECM. The OVs replicate inside cancer cells, leading to their lysis.}
    \label{cancer-ov-fibre}
\end{figure}

In this \re{study}, we develop further the non-local multiscale approach for modelling cancer-OV interactions introduced in \cite{Alsisi2020} \re{(and which focused on the role of cell-cell and cell-matrix adhesive interactions on the spread of OV throughout solid tumours) by combining it with the heterogeneous ECM approach proposed in} \cite{Shuttleworth2019}, to investigate these tumour-ECM-OV interactions in fibrous ECM; see Figure \ref{cancer-ov-fibre}.

In Section~\ref{Sect:Model3} we describe the new multiscale mathematical model. In Section~\ref{Sect:Numerics3} we describe the numerical approach used to discretise the macroscale and microscale equations, while in Section~\ref{Sect:Results3} we present the results of the numerical simulations. We conclude in Section~\ref{Sect:Summary3} with a summary and discussion of the results.   

\section{Mathematical Model}\label{Sect:Model3}
\begin{figure}[h]
    \centering
    \includegraphics[width=\textwidth]{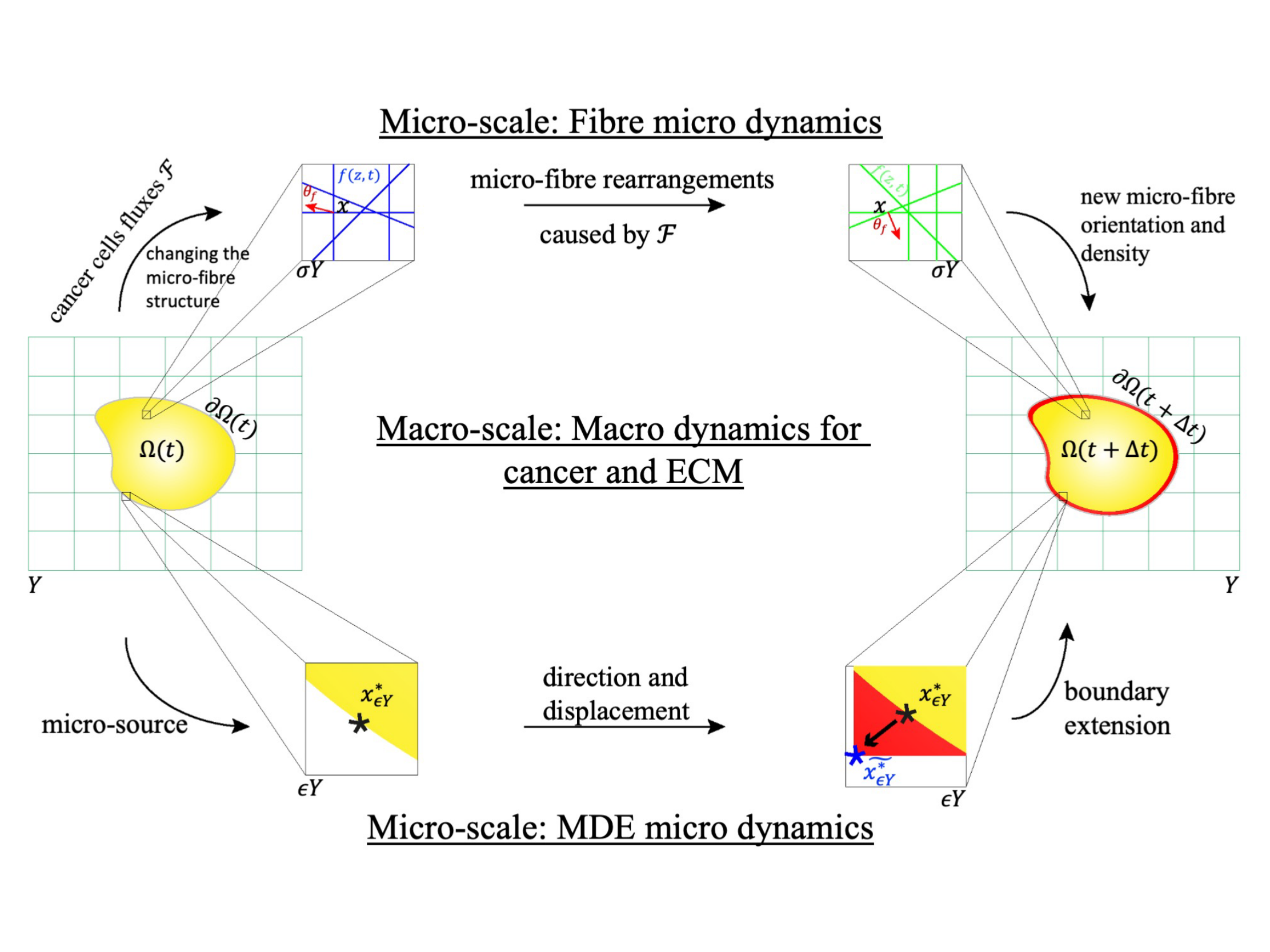}
    \caption{Schematic demonstrate the links between the multiscale model, specifically, the link between the macro-scale dynamics and the micro-fibre rearrangements, and the extension of cancer region $\Omega(t)$ caused by MDE micro dynamics.}
    \label{macro-micro-fibre}
\end{figure}

The multiscale moving boundary model used here is based on the two-scale (tissue scale -- macro-scale, and cell scale -- micro-scale) moving boundary framework \re{introduced in \cite{Trucu2013} and recently applied to nonlocal cell-cell interactions in the context of oncolytic viral therapies}~\cite{Alsisi2020}. Furthermore, in here we explore the dynamic interaction between an invading \re{solid} tumour, OV, and a two-component ECM (\re{which} was first introduced by \cite{Shuttleworth2019}). This complex dynamic is captured by two interconnected multiscale systems that share the same macro-scale cancer dynamics at the tissue-scale. However, at cell-scale we use two distinct micro-scale dynamics for fibre rearrangement and \re{for cancer invasion boundary, both being} linked to the macro-dynamics through two double feedback loops, as illustrated in Figure \ref{macro-micro-fibre}. 
\subsection{Macro-\dtr{S}cale Dynamics}
\re{Let us denote by} $Y \subset \R$ the macro-scale domain representing \re{the} maximal environmental tissue square \re{that we consider in this study}. \re{Also, let us denote by} $\Omega(t)\subset Y$ the spatial support of the growing tumour within the macro-scale domain, at a time $t\in[0,T]$. Let $c(x,t)$ and $i(x,t)$, $\forall (x,t)\in\Omega(t)\times[0,T],$ represent the spatio-temporal densities of the uninfected cancer cells and the infected cancer cells, respectively. \re{Finally, let} $v(x,t)$ and $e(x,t),$ $\forall (x,t)\in Y \times [0,T]$, represent the spatial-temporal densities of the oncolytic virus and the cumulative extracellular matrix, respectively. The later \re{is} defined as
\begin{equation*}
e(x,t) = E(x,t) + F(x,t),
\end{equation*} 
where $F(x,t)$ denotes the macroscale spatial distribution of the fibre ECM, which accounts for all significant ECM fibres such as collagen fibres or fibronectin fibrils, and $E(x,t)$ denotes the spatial distribution of the non-fibre ECM, which includes all the other non-fibre components of the ECM, i.e., elastin, laminins, fibroblasts, etc. At any spatio-temporal node $(x,t) \in Y \times [0, T]$, the macroscale dynamics for ECM \re{is} described by the following equations:
\begin{eqnarray}
\dfrac{\partial E}{\partial t} &=& -E(\alpha_{c}\,c+\alpha_i\, i) +\mu_2 E(1-\rho(\bf(s))),\label{eq:ecm}\\
\dfrac{\partial F}{\partial t} &=& -F(\alpha_{cF}\,c+\alpha_{iF} \,i),\label{eq:macro-fibre}
\end{eqnarray}
where $\alpha_{c},\alpha_{i}, \alpha_{cF}, \alpha_{iF}>0$ are the ECM degradation rates caused by cancer cells subpopulations $c$ and $i$, respectively. Further, $\mu_{2}>0$ is a fixed remodelling rate. The fibre micro-dynamics will be discussed in details in Section \ref{micro-fibre}. Defining $\textbf{s}$ the \emph{tumour-ECM vector,}
\begin{equation*}
\textbf{s}\left(x,t\right) = (c(x, t), i(x, t), e(x,t))^T,
\end{equation*}
then the volume fraction of space occupied by the tumour \re{and the ECM} is given as
\begin{equation}\label{fraction}
\rho (x,t) \equiv \rho(\textbf{s}(x,t)) := \nu_e e(x,t) + \nu_c \big(c(x,t)+i(x,t) \big),
\end{equation}
where $\nu_e$ represents the fraction of physical space occupied by the ECM and $\nu_c$ is the fraction of physical space occupied collectively by all cancer subpopulations.

In this study tumour dynamics consists of three main components: motility, proliferation, and death. \re{Moreover,} we describe \re{the spatial fluxes for} both cancer cells $(c,i)$ by a combination of a linear diffusion term caused by cells' random walk and a directed migration term due to cell-cell and cell-ECM adhesion. Here we investigate directed cell migration from two perspectives: local migration and non-local \re{migration as a result to} cell adhesion \cite{Bhagavathula2007,CAVALLARO200139,Humphries2006,Ko_2001,Wijnhoven2000,Chaplain2006}. \re{Therefore, the dynamics of} the uninfected cancer cells subpopulation is \re{given by}: 
\begin{equation}\label{c_eq1}
\frac{\partial c}{\partial t} =  \nabla \cdot \left[ D_c  \nabla c -c\,\A_c(t,x,\textbf{s},\theta_f)\right] + \mu_1 c(1-\rho(\textbf{s})) - \varrho cv,
\end{equation}
where $D_{c}>0$ is a constant diffusion coefficient, $\mu_{1}>0$ is a constant proliferation coefficient, $\varrho > 0$ is a constant rate at which the uninfected cancer population diminishes due to infection by the oncolytic virus $v$, while $\A_c(t,x,\textbf{s},\theta_f)$ is a non-local spatial flux term that describes the cells adhesion process that causes cancer cells to \re{move in a directed manner}. In \cite{Alsisi2020,aalsisi21} \re{the authors} studied the effects of cell adhesion process on cancer-OV interaction, by focusing on cell-cell adhesion and cell-ECM-non-fibres substrate. Here we adopt the modelling concept proposed by \cite{Shuttleworth2019}, to consider the essential role performed by the cell-fibres adhesive interaction. Denoting the cell-cell adhesion function by:
\begin{equation}\label{ad_fun1}
\T_c(x+y,t) = \textbf{S}_{cc} \,c(x+y,t) + \textbf{S}_{ci} \,i(x+y,t),
\end{equation}
where $\bS_{cc}>0$ and $\bS_ {ci}>0$ are the strengths of cell-cell adhesion and cross adhesion bonds, respectively  that are formed between cancer cells distributed at $x$ and cells at $x+y$. Since the cell-cell adhesion strength $\textbf{S}_{cc}, \textbf{S}_{ci}$ is dependent on the quantity of intercellular $Ca^{2+}$ ions available within the ECM \cite{Gu2014,Hofer2000}. We adopt a similar approach as in \cite{Shuttleworth2019} to compute cell-cell adhesion strength, and extend it to determine cross cell-cell adhesion strength as follows: 
\begin{equation}\label{S_cceq}
\textbf{S}_{\cdot \cdot}(E) = S_{\cdot \cdot}^{\,\textit{max}} \exp{\left( 1-\frac{1}{1-(1-E(x,t))^2}\right)}, 
\end{equation}
where $\textbf{S}_{\cdot \cdot} \in \lbrace \textbf{S}_{cc}, \textbf{S}_{ci}, \textbf{S}_{ic}, \textbf{S}_{ii} \rbrace$ with $S_{\cdot \cdot}^{\,\textit{max}} \in \lbrace S_{cc}^{\,\textit{max}}, S_{ci}^{\,\textit{max}}, S_{ic}^{\,\textit{max}}, S_{ii}^{\,\textit{max}}  \rbrace$ respectively. $S_{\cdot \cdot}^{\,\textit{max}}$ is fixed and represent the maximum strength of cell-cell adhesive junctions. Within a sensing radius $R>0$ at time $t>0$, the non-local adhesive flux is defined as follows:
\begin{IEEEeqnarray}{rCl}
\A_c(t,x,\textbf{s},\theta_f) & = &  \frac{1}{R} \int_{\textbf{B}(0,R)} \mathcal{K} (\| y\|_2)(1-\rho(\textbf{s}))^{+} \Big(\textbf{n}(y) [ \T_c(x+y,t) + \textbf{S}_{ce} \: E(x+y,t)]
\nonumber\\
&& +\>  \bar{\bn}(y,\theta_f(x+y,t)) \textbf{S}_{cF}\: F(x+y,t)\Big)\, \chi_{_{\Omega(t)}}(x+y,t)\,dy,
\label{adhesion_eq3}
\end{IEEEeqnarray}
where $\Bila(0, R):=\lbrace z\in \R^2:|z|\leq R\rbrace$ is a closed ball centred at origin and of radius R, called here the \emph{sensing region}. $\textbf{S}_{ce}$ is a constant \re{that} describes the cell-ECM adhesion strength. $\textbf{S}_{cF}$ is a constant \re{that} describes the cell-fibre-ECM adhesion strength. $\chi_{_{\Omega(t)}}(\cdot)$ is the characteristic function of $\Omega(t)$, and $(1-\rho(\textbf{s}))^{+} := \text{max} \lbrace (1-\rho(\textbf{s})),0\rbrace$ is a threshold term to avoid local overcrowding, and $\bn(y)$ denotes the unit radial vector giving by: 
\begin{equation}
\bn(y):= \begin{cases} 
      \dfrac{y}{\norm{y}_2} &  \text{if}\quad y \in \textbf{B}(0,R)\setminus \lbrace (0,0)\rbrace,\\
      (0,0) & otherwise,
   \end{cases}
\end{equation}
with $\nor{\cdot}_{_{2}}$ representing the usual Euclidean norm. $\bar{\bn}(y,\theta_f(x+y,t))$ is the unit vector that is dependent on the fibre orientations, defined as follows:
\begin{equation}
\bar{\bn}(y,\theta_f(x+y,t)):= \begin{cases} 
      \dfrac{y+\theta_f(x+y,t)}{\norm{y+\theta_f(x+y,t)}_2} &  \text{if}\quad y \in \textbf{B}(0,R)\setminus \lbrace (0,0)\rbrace,\\
      (0,0) & otherwise,
   \end{cases}
\end{equation}
where $\theta_f(x+y,t)$ is the orientation of the fibres at macro-scale, this was first introduced by \cite{Shuttleworth2019}. This orientation is derived by the micro-scale mass distribution of micro-fibres $f(.,t)$, in the sense that affects the cell-ECM adhesion and characterises the ECM fibres distributed at the macro-scale location $x\in Y$ namely $F(x,t)$; for detailed mathematical formulas see Section~\ref{micro-fibre}. \re{Note that} $F(x+y,t)$ \re{describes the} influence \re{of fibres distributed at $x + y$ on the adhesion of cells} at location $x$ (with adhesion strength $S_{cF}$); see Figure \ref{macro-micro-fibre} for more details. Furthermore, the radially symmetric kernel $\mathcal{K}(\cdot):[0,R]\to [0,1]$ explores the dependence of the strengths of the established cell adhesion junctions on the radial distance from the centre of the sensing region $x$ to $\zeta\in \Bila(x,R):= x + \Bila(0, R)$.  Since these adhesion junction strengths are assumed to decrease as the distance $r:=\nor{x-\zeta}_{_{2}}$ increases, $\mathcal{K}$ therefore is taken here of the form
\begin{equation}
\mathcal{K}(r):= \frac{3}{2\pi R^2}\left(1- \dfrac{r}{2R}\right).
\end{equation}
By summing up the radially distributed adhesive interactions between the cells at $x$, and the cells and ECM at $x+y$ within $y\in \Bila (x,R)$, the term $\frac{1}{R}$ that appears in the front of the expression (\ref{adhesion_eq3}) is simply an interaction range normalisation factor.

Next, \re{we focus on} the infected cancer cell subpopulation $i(x,t)$ that emerges within this dynamics due to infections by the OV. \re{We denote by} $\varphi_{i}(u)$ the effect of the cell adhesion processes \re{that take place} either locally (through adhesive interactions between infected cancer cells and ECM, as tumour cells exercise haptotactic movement towards higher levels of ECM), or non-locally (where both cell-cell and cell-ECM adhesive interactions are accounted for within an appropriate cell sensing region). Mathematically, this can be formalised as
\begin{equation}\label{new_f1}
\varphi_{i}(\bs) := \begin{cases} 
      \eta_{i} \!\nabla\!\! \cdot \!\big(i\nabla e\big), &   \footnotesize\text{local \dt{haptotactic interactions between infected cancer cells and ECM}}, \\
      \nabla \!\! \cdot\! \big( i\A_{i}(\cdot,\!\cdot,\! \bs(\cdot,\!\cdot)\!)\!\big), & \footnotesize\text{non-local \dt{cell$-$cells and cell$-$ECM interactions on a cell sensing region}},
   \end{cases}
\end{equation}
where $\eta_{i}$ is a constant haptotactic rate associated to $i$, while $\A_{i}\big(\cdot,\!\cdot,\! \bs(\cdot,\!\cdot) \!\big)$ is a non-local spatial flux term defined as in equation (\ref{adhesion_eq3}). \re{We define} $\T_i(x+y,t)$ as follows
\begin{equation}\label{ad_fun1}
\T_i(x+y,t) = \textbf{S}_{ic} \,c(x+y,t) + \textbf{S}_{ii} \,i(x+y,t),
\end{equation} 
where $\textbf{S}_{ic}$ and $\textbf{S}_{ii}$ are are the strengths of \re{infected} cell-cell adhesion and cross adhesion bonds, and are dependent on ECM as defined in equation (\ref{S_cceq}). Then, $\A_{i}\big(\cdot,\!\cdot,\! \bs(\cdot,\!\cdot)\!\big)$ becomes
\begin{IEEEeqnarray}{rCl}
\A_i(t,x,\textbf{s},\theta_f) & = &  \frac{1}{R} \int_{\textbf{B}(0,R)} \mathcal{K} (\| y\|_2)(1-\rho(\textbf{s}))^{+} \Big(\textbf{n}(y) [ \T_i(x+y,t) + \textbf{S}_{ie} \: E(x+y,t)]
\nonumber\\
&& +\>  \bar{\bn}(y,\theta_f(x+y,t)) \textbf{S}_{iF}\: F(x+y,t)\Big)\, \chi_{_{\Omega(t)}}(x+y,t)\,dy,
\label{adhesion_eq_i3}
\end{IEEEeqnarray}
where $\textbf{S}_{ie}$ is a constant describe the cell-ECM adhesion strength and $\textbf{S}_{iF}$ is a constant describe the cell-fibre-ECM adhesion strength. Thus, the \re{spatio-temporal dynamics of the} infected cancer cell subpopulation is governed by the following equation
\begin{equation}\label{i_m_equ1}
\dfrac{\partial i}{\partial t} = D_{i} \Delta i-\dt{\varphi_{i}(\bs)} + \varrho c v -\delta_i i.
\end{equation}
Here $D_{i}>0$ is a constant random motility coefficient, and $\varphi_{i}(\bs)$ represents the directed migration induced by the cell-adhesion processes that corresponds to $i$ and is described in \eqref{new_f1} Further, the infected population expands at a rate $\varrho$ due to new infections occurring among the uninfected cells, and dies at rate $\delta_{i}>0$.

Next, for the oncolytic virus spatio-temporal dynamics, we adopt here a similar reasoning as in \cite{Alsisi2020,aalsisi21}, \re{and assume} that the OV motion is described by a random movement that is biased by a "haptotactic-like" spatial transport towards higher ECM levels. Thus, the dynamics of the oncolytic virus is governed by
\begin{equation} \label{v_equ1}
\dfrac{\partial v}{\partial t} = D_v \Delta v -\eta_v \nabla \cdot \left(v\nabla e\right) + b i - \varrho c v - \delta_v v.
\end{equation}
Here $D_{v}>0$ is a constant random motility coefficient, $\eta_{v}>0$ is a constant haptotactic coefficient, $b>0$ is a viral replication rate within infected cancer cells, and $\delta_{v}>0$ is the viral death rate. 

Finally, the coupled interacting tumour-OV macro-dynamics is governed by \eqref{c_eq1}-\eqref{v_equ1}
in the presence of initial conditions
\bequ
c(x,0)=c^{0}(x), \quad i(x,0)=i^{0}(x), \quad  v(x,0)=v^{0}(x), \quad \forall x \in \Omega(0), 
\eequ
while assuming zero-flux boundary conditions at the moving tumour interface $\partial \Omega(t)$.

\subsection{Fibre Micro Dynamics \dtr{on the Bulk of the Tumour}}\label{micro-fibre}
The cancer cells macro-dynamics cause the fibres to undergo a microscopic rearrangement process in addition to the macroscale fibre degradation mentioned in equation (\ref{eq:macro-fibre}). We start by defining the micro-fibre domain as $\sigma Y(x) := x + \sigma Y, \forall x \in Y$ with scale size $\sigma>0$ following \cite{Shuttleworth2019}, in which the fibre ECM rearrangement occurs reflecting on fibre macroscale orientation as illustrated in Figure \ref{macro-micro-fibre}. At any macroscale point $x \in \Omega(t)$, the ECM-fibre phase is described by a macroscale vector field $\theta_{f}(x, t)$, which give us the amount of fibres distributed at $(x, t)$ and their macroscopic fibres orientation computed by the revolving barycentral orientation $\theta_{f, \sigma Y(x)}(x, t)$ generated by the microscopic mass distribution of microfibres $f(\cdot, t)$ within the micro-domain $\sigma Y(x)$. Based on the approach introduced by \cite{Shuttleworth2019}, the macroscale fibre orientation $\theta_{f}(x, t)$ is defined as follows:  
\begin{equation}
\theta_{f}(x, t)=\frac{1}{\lambda(\sigma Y(x))} \int_{\sigma Y(x)} f(z, t) d z \cdot \frac{\theta_{f, \sigma Y(x)}(x, t)}{\left\|\theta_{f, \sigma Y(x)}(x, t)\right\|_{2}},
\end{equation}
where $\lambda(\cdot)$ is the usual Lebesgue measure, $\|\cdot\|_{2}$ represents the usual Euclidean norm, and 
\begin{equation}
\theta_{f, \sigma Y(x)}(x, t)=\frac{\int_{\sigma Y(x)} f(z, t)(z-x) d z}{\int_{\sigma Y(x)} f(z, t) d z},
\end{equation}
is the naturally generated revolving barycentral orientation $\theta_{f, \delta Y(x)}(x, t)$ associated with $\sigma Y(x)$ is given by the Bochner-mean-value of the position vectors function $\sigma Y(x) \ni z \mapsto z-x \in \R^{N}$ with respect to the density measure $f(x, t) \lambda(\cdot)$. The macroscopic mean-value fibre representation at any $(x, t)$ is then given by the Euclidean magnitude of $\theta_{f}(x, t)$, namely,
\begin{equation}
F(x, t):=\left\|\theta_{f}(x, t)\right\|_{2}.
\end{equation}

At any time $t$ and at any spatial location $x \in Y$, the cancer cells realign the micro-fibres through a microscopic rearrangement process in each micro-domain $\sigma Y(x)$ that is triggered by the combined macroscale spatial flux of both cancer cell subpopulations as follows:
$$
\mathcal{F}(x, t):=\mathcal{F}_{c}(x, t)+\mathcal{F}_{i}(x, t)
$$
where
\begin{eqnarray}
\mathcal{F}_{c}(x, t)&:=&D_c \nabla c(x, t)-c(x, t) \mathcal{A}_{c}\left(x, t, \bs(\cdot, t), \theta_{f}(\cdot, t)\right) \\
\mathcal{F}_{i}(x, t)&:=&D_{i} \nabla i(x, t)- \varphi_{i}(\bs)\label{flux_i},
\end{eqnarray}
where 
\begin{equation}\label{new_f2}
\bar{\varphi_{i}}(\bs) := \begin{cases} 
      i\nabla e, \\
       i\A_{i}(x,t,\bs(\cdot,t),\theta_{f}(\cdot, t)).
   \end{cases}
\end{equation}
For simplicity, denoting the \emph{total cancer cell population} by $c_{\textrm{total}}(x,t) = c(x,t)+i(x,t)$. The combined flux $\mathcal{F}(x, t)$ acts upon the micro-scale distribution $f(z, t), \forall z \in \sigma Y(x)$ in accordance to the magnitude that the total mass of cancer cells has relative to the combined mass of cells and fibres at $(x, t)$, which is given by the weight
\begin{equation}
\omega(x, t)=\frac{c_{\textrm{total}}(x,t)}{c_{\textrm{total}}(x,t)+F(x, t)}.
\end{equation}

\begin{figure}[h]
    \centering
    \includegraphics[width=\textwidth]{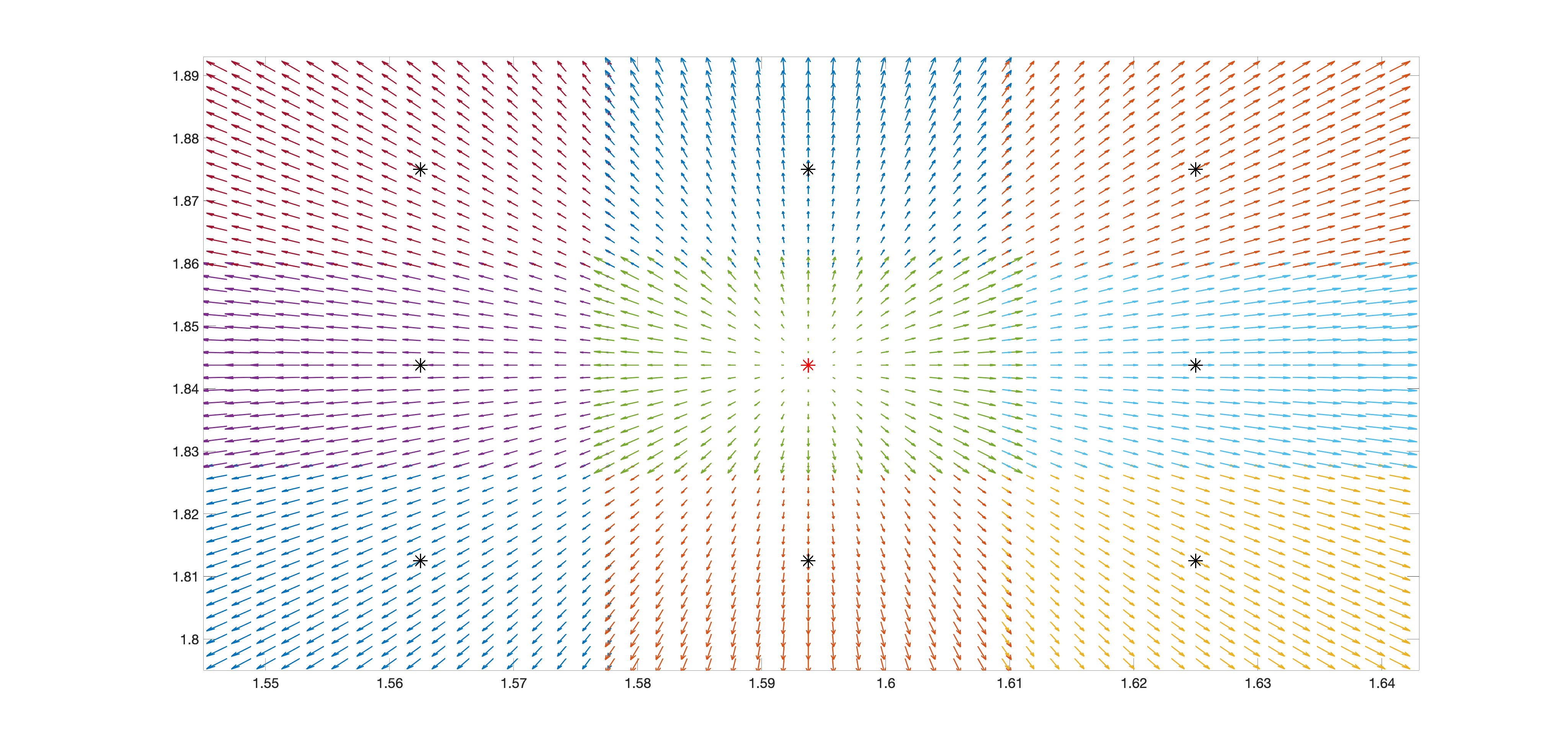}
    \caption{Illustration to show $\sigma Y(x)$ with $x$ is red star and the green vector field shows $x_{\mathrm{dir}}(z)$. The other black stars are the neighbours of point $x$ with their respected $\sigma Y$ to show full direction from the centred $x$ to all $z$ belong to $\sigma Y(x)$ and its surrounding boxes.}
    \label{xdir}
\end{figure}

At the same time, the \dtr{total} spatial flux of cancer cells $\mathcal{F}(x, t)$ is balanced in a weighted manner by the orientation $\theta_{f}(x, t)$ of the existing distribution of fibres at $(x, t)$ that is appropriately magnified by a weight that accounts for the magnitude of fibres versus the combine mass of cells and fibres at $(x, t)$ and is given by $(1-\omega(x, t))$. As a consequence, the micro-scale distribution of micro-fibres $f(z, t), \forall z \in \sigma Y(x)$ is therefore acted upon uniformly by the resultant force given by the following rearrangement  vector-valued function
\begin{equation}
r(\sigma Y(x), t):=\omega(x, t) \mathcal{F}(x, t)+(1-\omega(x, t)) \theta_{f}(x, t).
\end{equation}
\dtr{In this context, a} mass distribution of the micro-fibres $f(z,t)$ on $\sigma Y(x)$ is \dtr{exercised under the influence of} this rearrangement vector $r(\sigma Y(x), t)$, \dtr{resulting in spatial relocation of micro-fibres on both} $\sigma Y(x)$ and \dtr{its} neighbouring micro-domains. \dtr{Indeed, under the incidence of $r(\sigma Y(x), t)$, a certain fraction of the micro-fibres positioned at a given $z\in \sigma Y(x)$ get transported at} new micro-scale position $z^{*}$, \dtr{given by}
\begin{equation}\label{new_fibre_pt}
z^{*}:=z+\nu_{\sigma Y(x)}(z, t),
\end{equation}
where $\nu_{\sigma Y(x)}(z, t)$ is the \dtr{emerging} relocation vector: 
\begin{equation}\label{relocation_vector}
\nu_{\sigma Y(x)}(z, t)=\left(x_{\mathrm{dir}}(z)+r(\sigma Y(x), t)\right) \cdot \frac{f(z, t)\left(f_{\max }-f(z, t)\right)}{f^{*}\dtr{(z,t)}+\left\|r(\sigma Y(x))-x_{\mathrm{dir}}(z)\right\|_{2}} \cdot \chi_{\{f(\cdot, t)>0\}}(z).
\end{equation}
\dtr{Here, $x_{\mathrm{dir}}(z)=\overrightarrow{xz}$ is the barycentric position vector pointing to $z$ in $\sigma Y(x)$ see Figure \ref{xdir}, $f_{max}$ represents the maximum level of fibres that could reside at the micro-location $z\in \sigma Y$ at any given time, and $f^{*}:=f(z,t)/f_{max}$ is the local micro-fibres saturation fraction. Finally, the relocation magnitude in the direction $x_{\mathrm{dir}}(z)+r(\sigma Y(x), t)$ is simultaneously mediated by the ability of micro-fibres to dislocate (which is exercised when these are not at the maximum level) and is attenuated by the level of the micro-fibres mass fraction $f^{*}$ at $z$ in conjunction with the barycentric position defect quantified here by $\left\|r(\sigma Y(x))-x_{\mathrm{dir}}(z)\right\|_{2}$.} 
Therefore, provided that micro-fibres levels at $z$ are not at their maximum level $f_{\max }$ a micro-fibres mass transport from $z$ to the location $z^*$ is exercised, while lower levels of micro-fibres saturations at $z$ combined with a better \dtr{flux} alignment provided by a smaller position defect result in a relocation of the micro-fibres mass in direction $(x_{\mathrm{dir}}(z) + r(\sigma Y(x), t))$ at a greater distance. Finally, this transport \dtr{occur provided that a  enough room is available at position $z$. This being is naturally mediated through the} movement probability
\begin{equation}
P_{\textrm{move}} := \max \left( 0, \frac{f_{\max} - f(z^*,t)}{f_{\max }} \right),
\end{equation} 
that \dtr{quantifies the capacity that is still available at $z^{*}$, and} enables only the amount $f(z,t)P_{\textrm{move}}$ of micro-fibres to be transported to the new location and the rest $f(z,t)(1-P_{\textrm{move}})$ remain at their location. 

\subsection{\dtr{MDEs} Microscale Dynamics \dtr{at the Leading Edge and the Induced Tumour Boundary Movement}}
\dtr{In addition to the fibres micro-dynamics that occures in the bulk of the tumour, a second type of tumour micro-dynamics is generated by the cell scale activity of the matrix degrading enzymes (MDEs) along the tumour invasive edge \cite{Hanahan_et_al_2011,Weinberg2006}. Indeed, secreted by the cancer cell population from the outer proliferating rim of the tumour, the MDEs are transported across the tumour interface within a cell-scale (micro-scale) neighbourhood of the tumour boundary in the peritumoural region, where they cause degradation to the ECM distribution that they meet, giving rise to further tumour progression. To capture this process mathematically, we adopt here a similar approach to the one developed in  \cite{Trucu2013}, whereby the entire cross-interface MDEs transport micro-process along $\partial \Omega(t)$ is decomposed into a union of boundary micro-processes that take place on a appropriately constructed covering bundle $\{\epsilon Y\}_{\epsilon Y\in \P_{\epsilon}(t)}$ of $\epsilon-$size overlapping micro-domains $\epsilon Y$, whose union form the cell scale neighbourhood $\Ne_{\epsilon}(\partial \Omega(t))$, \emph{i.e.}, 
\[
\Ne_{\epsilon}(\partial \Omega(t)):=\bigcup\limits_{\epsilon Y\in \P(t)} \epsilon Y
\]
where the entire enzymatic process and its consequences are explored. Thus, to explore the MDEs enzymatic activity over a small time range $\Delta t>0$, denoting by $m(y,\tau)$ the micro-scale density of MDEs at the micro-scale spatio-temporal location $(y,\tau)\in \epsilon Y \times [0, \Delta t]$, we have that a micro-scale source of MDEs $ \mathcal{G}_{\epsilon Y}:\epsilon Y \times [0, \Delta t] \to [0,\infty)$ is induced naturally as a collective contribution of the cancer cells that are within a small distance $\rho>0$ from each $y\in \epsilon Y$, and so this can be formalised mathematically as 
\begin{equation}\label{source_3}
 \mathcal{G}_{\epsilon Y} (\dt{z},\tau) =
 \begin{cases}
\dfrac{\int\limits_{\textbf{B}(z,\rho)\cap\Omega (t_0)} \gamma_{c}c (x, t_{0}+\tau)+\gamma_{i}i (x, t_{0}+\tau)dx} {\lambda (\textbf{B}(z,\rho)\cap\Omega (t_0))} , & z\in \epsilon Y \cap \Omega (t_0),\\
 0, & otherwise,
 \end{cases}
\end{equation}
Here $\lambda(\cdot)$ is the standard Lebesgue measure on $\mathbb{R}^N$, and $\Bila(z,r):= \lbrace x \in Y: \norm{z-x}_\infty \leq \rho\rbrace$ is the small active tumour region closed to the tumour interface where the cancer cells collectively contribute to the creation of the source of MDEs at $z\in \epsilon Y$ over the time interval $([t_{0},t_{0}+\Delta t]$. Finally, and $\gamma_{c}>0$, $\gamma_{i}>0$ represent constant MDEs secretion contributions of the uninfected and infected cancer cells, respectively. Finally, in the presence of the MDEs source  $\mathcal{G}_{\epsilon Y} (\dt{z},\tau)$, the MDES are assumed here to exercise a diffusive transport within the entire micro-domain $\epsilon Y$, which is mathematically formulated through the following reaction-diffusion equation }
\begin{equation}\label{MDE_equ3}
\dfrac{\partial m}{\partial \tau} = D_m \Delta m+ \mathcal{G}_{\epsilon Y} (z,\tau),  \quad\quad z\in \epsilon Y,\,\tau \in [0,\Delta t].
\end{equation}
\dtr{Furthermore, as no memory of pre-existing distributions of MDEs are assumed for the enzymatic process, and no molecular transfer is assumed across the boundaries of $\epsilon Y$, the proteolytic boundary micro-dynamics \eqref{MDE_equ3} that takes place on each $\epsilon Y$ is assumed to take place with null initial conditions and  flux-zero boundary conditions, i.e., 
\begin{equation}\label{IC_BC_microDynamics3}
\begin{array}{lll}
m(z,0) &=& 0, \\[0.2cm]
\bn_{_{\epsilon Y}}\cdot \nabla m\bigg |_{\partial \epsilon Y} &=&0,
\end{array}
\end{equation}
where $\bn_{_{\epsilon Y}}$ is the usual outward unit normal direction on $\partial\epsilon Y$.}

\dtr{Finally, the pattern of ECM degradation within the peritumoural region of the microscale neighbourhood $Ne_{\epsilon}(\partial \Omega(t))$ will correspond to the pattern of significant MDEs transport within each micro-domain $\epsilon Y$. In this context, following the multiscale mathematical modelling approach developed in  \cite{Trucu2013}, we obtain the law for the macro-scale boundary movement, which is specified in terms of direction and displacement magnitude for the spatial relocation of the tumour interface $\epsilon \cap \Omega Y$ for each boundary micro-domain. }
\subsection{Summary of Multiscale Model} 
In summary, the multiscale moving boundary model that we obtained for the tumour -- OV --ECM interactions  (see Figure \ref{macro-micro-fibre} \dtr{for an illustrative sketch}) \dtr{consists of the following parts:}
\begin{subequations}\label{macro_micro_s3}
\begin{align}
\textrm{the tumour $-$ OV $-$ ECM}&\;\; \textrm{macro-dynamics:} \nonumber\\[-0.2cm]
\frac{\partial c}{\partial t} &=  \nabla \cdot \left[ D_c  \nabla c -c\,\A_c(t,x,\textbf{s},\theta_f)\right] + \mu_1 c(1-\rho(\textbf{s})) - \varrho cv,\label{c_eq3_s3}\\
\dfrac{\partial i}{\partial t} &= D_{i} \Delta i-\dt{\varphi_{i}(\bs)} + \varrho c v -\delta_i i,\label{i_eq3_s3}\\
\dfrac{\partial E}{\partial t} &= -E(\alpha_{c}\,c+\alpha_i\, i) +\mu_2 E(1-\rho(\bf(s))),\label{eq:ecm_s3}\\
\dfrac{\partial F}{\partial t} &= -F(\alpha_{cF}\,c+\alpha_{iF} \,i),\label{eq:macro-fibre_s3}\\
\dfrac{\partial v}{\partial t} &= D_v \Delta v -\eta_v \nabla \cdot \left(v\nabla e\right) + b i - \varrho c v - \delta_v v,\label{virus:eq3_s3}\\[0.5cm]
\textrm{\dtr{micro-dynamics of fibres}}&\;\; \textrm{\dtr{triggered by the macro-scale-induced rearrangement flux}}\\[0.5cm]
\dtr{r(\sigma Y(x), t)}&\dtr{:=\omega(x, t) \mathcal{F}(x, t)+(1-\omega(x, t)) \theta_{f}(x, t)}\\[0.5cm]
\textrm{\dtr{induced on each $\sigma Y$ by}}&\;\; \textrm{\dtr{the total cell flux $\mathcal{F}(x, t)$, with $\omega(x, t):=\frac{c_{\textrm{total}}(x,t)}{c_{\textrm{total}}(x,t)+F(x, t)}$}}\\[0.5cm]
\textrm{the MDEs boundary} \textrm{micro}&\textrm{-dynamics:} \nonumber\\[-0.2cm]
\dfrac{\partial m}{\partial \tau} &= D_m \Delta m+ \mathcal{G}_{\epsilon Y} (z,\tau)\label{MDE_equ_s3}
\end{align}
\end{subequations}
The macro-dynamics and \dtr{the two} micro-dynamics are connected through \dtr{two} double feedback loop\dtr{s}:
\begin{itemize}
\item a \dtr{tumour bulk} \emph{top-down} link by which the macro-scale \dtr{cell flux triggers the micro-scale fibres rearrangement} 
\item a dtr{\emph{bottom-up} fibres link by which the micro-scale distribution of fibres naturally induces a spatial orientation that alters the tumour macro-dynamics}
\item a \dtr{leading edge \emph{top-down} link through which the tumour macro-dynamics induces a micro-scale MDEs molecular source in the cell-scale neighbourhood of $\partial \Omega(t)$}
\item a \emph{bottom-up} link \dtr{through} which the MDEs micro-dynamics \dtr{induces the law for}  macro-scale tumour boundary movement. 
\end{itemize}

\section{\dtr{Brief Overview of the} Computational Approach \dtr{for a Couple of Tumour-OV-ECM Interaction Scenarios}}\label{Sect:Numerics3} 
The numerical approach and computational implementation of the novel multiscale moving boundary model \dtr{proposed in this work builds directly} on the multiscale moving boundary computational framework initially introduced by \cite{Trucu2013} and further expanded in \cite{Alsisi2020,aalsisi21,Shuttleworth2019}\dtr{, and so while for full details we defer the reader to these references, in the following we will give a brief summary of this approach.}

\paragraph{Macro-scale numerical approach}
\dtr{For all the computational simulations we consider t}he macro-scale \dtr{tissue} domain $Y [0,4]\times [0,4]$, \dtr{and we discretise this} \dtr{with a uniform grid} $Y^{d}:=\{(x^{1}_{i},x^{2}_{j})\}_{i,j=1\dots N}$, with $N=[4/h]+1$, \dtr{of spatial step size $\Delta x = \Delta y:=h$, with $h>0$}. \dtr{Correspondingly, the discretised version of the tumour support is denoted by}  $\Omega^{d}(t)$ (i.e., $\Omega^{d}(t)=Y^{d}\cap \Omega(t)$), \dtr{with} $\partial \Omega^{d}(t)$ \dtr{standing for the tumour boundary}.
\dtr{In brief, for the macro-scale part of the code, we used a method of lines type approach, with the non-local time marching scheme developed in  \cite{Shuttleworth2019}. Finally, as detailed in \cite{Trucu2013}, the discrete domain $\Omega^{d}(t)$ of the progressing tumour is appropriately evolved with additional spatial nodes corresponding to the new locations reached by the invading cancer}.

\paragraph{Approximating the two micro-dynamic processes that occur simoultaneously}
\dtr{In this model we have two types of micro-scale processes that link to the same macro-dynamics, namely the MDEs boundary micro-dynamics at the leading edge of the tumour and fibre micro-dynamics. Thus, for the MDEs boundary micro-dynamics part, proceeding as detailed in \cite{Trucu2013}, the boundary MDEs micro-dynamics is approximated involving central finite differences for the spatial operators on each micro-domain $\epsilon Y$, with a backward Euler time marching scheme. Finally, for the fibres micro-dynamics and their rearrangements, the implementation follows closely the modelling details given in Section \ref{micro-fibre}. Thus, for that part, considering the fibres micro-scale taking place on square micro-domains $\sigma Y$ that have their vertices at the dual mesh nodes, we proceed with the evaluation of the spatial flux operators and rearrangement vectors as well as the final numerical inference of the emerging spatial fibres orientation, as detailed in \cite{Shuttleworth2019}.} 
\subsection{Local vs Non-Local Directed Migration Due to Cell Adhesion}\label{localvsnonlocal3}
Our simulations of the multiscale model explore two distinct \re{scenarios regarding the directed migration of infected cancer cell subpopulation:} local and non-local migration at macro-scale (see eq.~\eqref{i_eq3_s3}). Specifically, we consider the following cases:     
\begin{enumerate}
\item \emph{Local advective flux for the infected cancer cells:} in \eqref{i_eq3_s3} we have $\varphi_{i}(\bs)=\eta_{i} \!\nabla\!\! \cdot \!\big(i\nabla e\big)$, and thus the macro-dynamics takes the following form:
\begin{subequations}\label{local_sys_3}
\begin{align}
\frac{\partial c}{\partial t} &=  \nabla \cdot \left[ D_c  \nabla c -c\,\A_c(t,x,\textbf{s},\theta_f)\right] + \mu_1 c(1-\rho(\textbf{s})) - \varrho cv,\\
\dfrac{\partial i}{\partial t} &= D_{i} \Delta i-\eta_{i} \!\nabla\!\! \cdot \!\big(i\nabla e\big) + \varrho c v -\delta_i i,\\
\dfrac{\partial E}{\partial t} &= -E(\alpha_{c}\,c+\alpha_i\, i) +\mu_2 E(1-\rho(\bf(s))),\\
\dfrac{\partial F}{\partial t} &= -F(\alpha_{cF}\,c+\alpha_{iF} \,i),\\
\dfrac{\partial v}{\partial t} &= D_v \Delta v -\eta_v \nabla \cdot \left(v\nabla e\right) + b i - \varrho c v - \delta_v v.
\end{align}
\end{subequations}
\re{Since $\bar{\varphi_{i}}(\bs) =  i\nabla e,$ the} macroscale spatial flux for $\mathcal{F}_i$ in equation (\ref{flux_i}) becomes
\begin{equation}
\mathcal{F}_{i}(x, t):= D_{i} \nabla i(x, t)- i\nabla e.
\end{equation}
\re{The numerical} results for this case \re{are} shown in Figures \ref{baseline_res3}-\ref{high_fibre_ad}. 
\vspace{0.4cm}
\item \emph{Nonlocal advective flux for the infected cancer cells: in} \eqref{i_eq3_s3} we have $\varphi_{i}(\bu)= \nabla \!\! \cdot\! \big( i\A_{i}(\cdot,\!\cdot,\! \bs(\cdot,\!\cdot)\!)\!\big)$, and thus the macro-dynamics take the following form:
\begin{subequations}\label{nonlocal_sys_3}
\begin{align}
\frac{\partial c}{\partial t} &=  \nabla \cdot \left[ D_c  \nabla c -c\,\A_c(t,x,\textbf{s},\theta_f)\right] + \mu_1 c(1-\rho(\textbf{s})) - \varrho cv,\\
\dfrac{\partial i}{\partial t} &=  \nabla \cdot \left[ D_i  \nabla i -i\,\A_i(t,x,\textbf{s},\theta_f)\right] + \varrho c v -\delta_i i,\\
\dfrac{\partial E}{\partial t} &= -E(\alpha_{c}\,c+\alpha_i\, i) +\mu_2 E(1-\rho(\bf(s))),\\
\dfrac{\partial F}{\partial t} &= -F(\alpha_{cF}\,c+\alpha_{iF} \,i),\\
\dfrac{\partial v}{\partial t} &= D_v \Delta v -\eta_v \nabla \cdot \left(v\nabla e\right) + b i - \varrho c v - \delta_v v,
\end{align}
\end{subequations}
\re{Since $\bar{\varphi_{i}}(\bs) =  i\A_{i}(x,t,\bs(\cdot,t),\theta_{f}(\cdot, t))$, the} macroscale spatial flux for $\mathcal{F}_i$ in equation (\ref{flux_i}) becomes
\begin{equation}
\mathcal{F}_{i}(x, t):= D_{i} \nabla i(x, t)- i\A_{i}(x,t,\bs(\cdot,t),\theta_{f}(\cdot, t)),
\end{equation}
\re{The numerical} results for this case \re{are} shown in Figures \ref{2nonlocal_fibre20}-\ref{2nonlocal_cross_ad}.
\end{enumerate}

\subsection{Initial Conditions}
\begin{figure}
    \centering
    \includegraphics[width=0.8\textwidth]{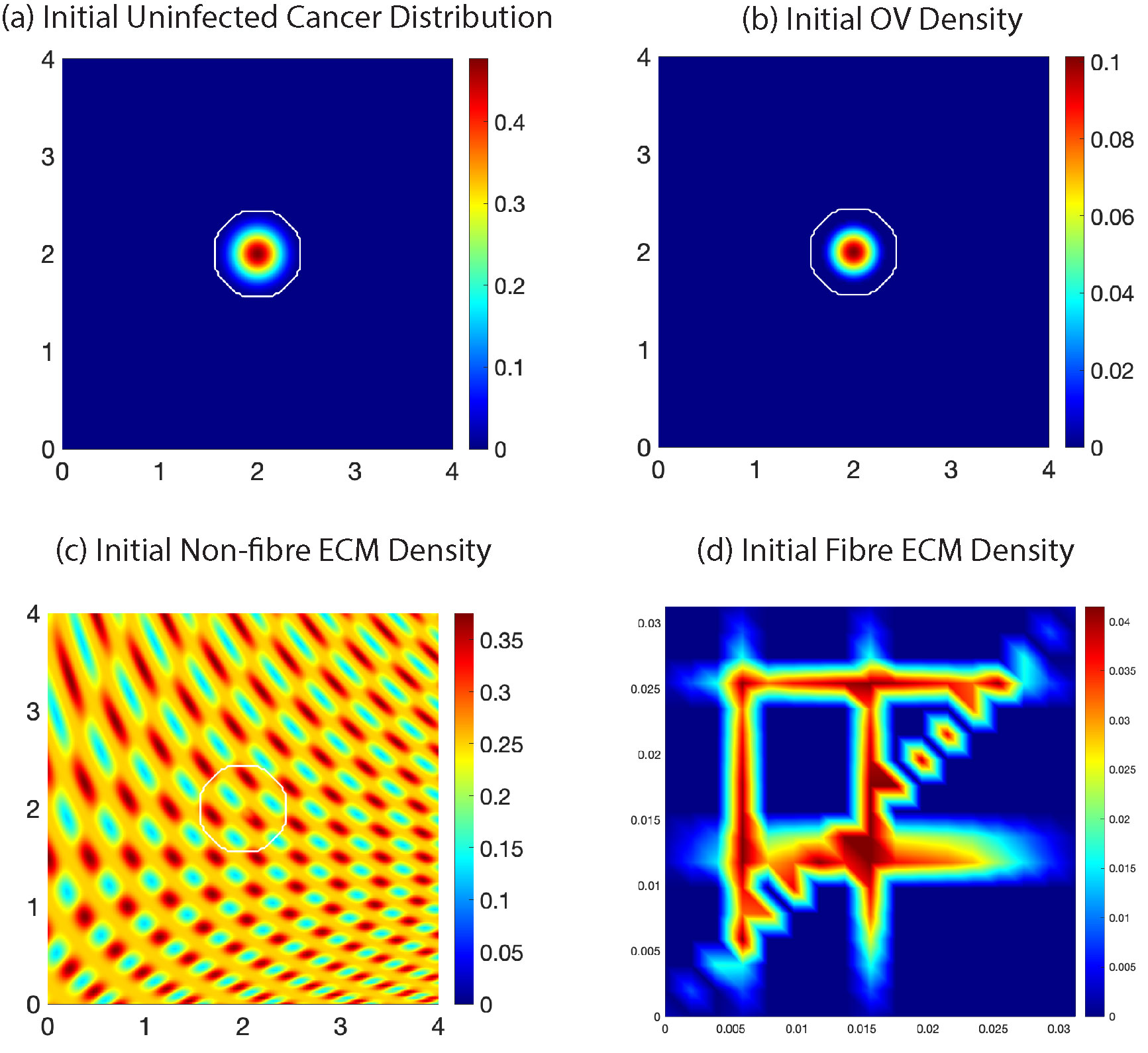}
    \caption{Initial conditions used for the numerical simulations: (a) uninfected cancer cells density ($c$), as given by equation (\ref{IC_c3}); (b) OV density ($v$), as given by equation (\ref{IC_v3}); (c) non-fibre ECM density ($E$), as given by equation (\ref{IC_e3}); (d) one micro-fibre domain which is repeated for every point on the macro-scale. The white curve \re{in sub-panels (a)-(c)} indicates the tumour boundary.}
    \label{IC3}
\end{figure}

The initial condition for the uninfected cancer cell population, $c(x,0)$ is chosen to describe a small localised pre-existing tumour aggregation. This is given by 
\begin{equation}\label{IC_c3}
c^{0}(x) = 0.5 \left(\exp\left(-\dfrac{\norm{x-(2,2)}_2^2}{2h}\right)-\exp\left(-3.0625 \right)\right)\left(\chi_{_{\mathbf{B}((2,2),0.5-\gamma)}}*\psi_\gamma\right), \quad \dt{\forall\, \, x\in Y,}
\end{equation}
whose plot is shown in Figure \ref{IC3}(a). Here $\psi_{\gamma}:\mathbb{R}^N \rightarrow \mathbb{R}_{+}$ is \dt{the usual} standard mollifier of radius $\gamma << \frac{\Delta x}{3}$ \dt{given by}
\begin{equation}
\psi_\gamma (x) :=  \dfrac{1}{\gamma^N} \psi\left(\dfrac{x}{\gamma}\right),  
\end{equation}
\re{with} the smooth compact support function $\psi$ given by
\begin{equation}
\psi(x):= \begin{cases} 
      \exp\frac{1}{\norm{x}_2^2-1} &  \text{if}\quad \norm{x}_2 < 1, \\
      0 & \text{otherwise}.
   \end{cases}
\end{equation}
We assume no infection at this stage, i.e. zero infected cancer cells ($i(x,0)$):
\bequ\label{IC_i}
i^{0}(x)=0, \;\;\;  \;\; \forall \; x\in Y.
\eequ
The initial condition for the non-fibre ECM density, $E(x,0)$, is given by an arbitrarily chosen heterogeneous pattern described by the following equations (as in~\cite{Shuttleworth2019})  
\begin{equation}
\label{IC_e3}
E(x,0) = \frac{1}{2} \,min\lbrace h(\zeta_1(x), \zeta_2(x)), 1-\dt{c^{0}_{p}(x)}\rbrace,
\end{equation}     
and is shown in Figure \ref{IC3}(c). Here, we have
\bequ
\begin{array}{rll}
\dt{h(\zeta_1(x), \zeta_2(x))} &:=& \frac{1}{2} +\frac{1}{4} \sin (\xi \dt{\zeta_1(x) \zeta_2(x)})^3 \cdot \sin \left(\xi \dt{\dfrac{\zeta_2(x)}{\zeta_1(x)}} \right),\\
\dt{(\zeta_1(x), \zeta_2(x))} &:=& \frac{1}{3} (x+\frac{3}{2}) \in [0,1]^2,  \quad \forall x\in Y, \quad and \quad \xi = 7\pi.
\end{array}
\eequ
\re{Furthermore,} the initial condition for one micro-scale fibre domain $\sigma Y(x)$ is shown in Figure \ref{IC3}(d), and it is repeated for all macro-scale locations. \re{To determine the density of the fibres at any macro-point $x$, we integrate the corresponding fibre-micro domain $\sigma Y(x)$.} Due to visibility reason, we avoid presenting the pattern of the microscale fibres on the macroscale.\re{ For the baseline simulations presented here, we choose the ratio of fibres to non-fibres components of ECM at $20\%:80\%$. }

\re{Finally,} the initial condition for the OV population ($v(x,0)$) is chosen to describe one single injection in the middle of the tumour aggregation, as in~\cite{Alsisi2020,aalsisi21}:
\begin{equation}\label{IC_v3}
v^{0}(x)= \Phi(x)\cdot \theta(v), 
\end{equation}
where
\bequ
\begin{array}{lll}
\Phi(x) &=& \frac{1}{8}\left(\exp\left(-\dfrac{\norm{x-(2,2)}_2^2}{2h}\right)-\exp\left(-1.6625 \right)\right), \\
\textrm{and}&&\\
\theta(v)&=&
\begin{cases} 
      1 &  \text{if}\quad \Phi(x) > 5\times10^{-5}, \\
      0 & \text{otherwise},
   \end{cases}\label{theta}
\end{array}
\eequ
\dtr{which is smoothed out on the frontier of the viral density support  $\Gamma_{v}:=\partial \{x\in Y\,|\, v^{0}(x)>0\}$.}


\section{Results}\label{Sect:Results3}
The numerical results shown in this section are computed using the parameter values listed in Table \ref{table_test3}, which we refer to as 'baseline parameters' for ease of reference. \re{Whenever we change these parameters,} we clearly specify the new values we use for the simulations. 

We start in Section \ref{res:baseline3} by investigating numerically the impact of fibre-ECM local approach for the infected cancer cells $i$ used to describe the cell-cell and cell-matrix adhesion flux on cancer-OV interaction. Then, in Section \ref{res:fibre_amount}, we investigate the impact of varying the amount of fibres in the ECM ($\mathcal{R}_F$). Next, in Section \ref{res:adhesion_strength3} we increase the cell-fibre-ECM adhesion strength \re{for} some of the \re{cases studied in the} previous sections. In Section \ref{res:2nonlocal} we investigate numerically the impact of fibre-ECM non-local approach for the infected cancer cells $i$ used to describe the cell-cell and cell-matrix adhesion flux on cancer-OV interaction and comparing it to a variety of distinct adhesion strengths. \re{In} section \ref{res:amount_fibre_2nonloca} \re{we investigate the impact} of increasing the amount of fibre in ECM \re{for} the non-local system. Finally, in Section \ref{res:cross_adhesion} we investigate the impact of different cross cell-cell adhesion strengths. 
 
\begin{table}[!h]
\caption{\re{Nondimensional} baseline parameters values used for our multiscale computations. \re{These baseline parameters were obtained from other articles (see references in the last column) or our own estimates.} }
\label{table_test3}
\renewcommand\arraystretch{1}
\noindent\[
\begin{array}{|c|c|l|l|}
\hline
\text{Param.} & \text{Value} & \text{Description} & \text{Reference}\\
\hline
D_{c}   & 0.00035    & \text{Uninfected cancer cell diffusion coefficient} & \text{\cite{Domschke2014a}}   \\
D_i    & 0.0054      & \text{Infected  cancer cell diffusion coefficient} &   \text{\cite{camara2013}}    \\
D_v  &  0.0036     & \text{Constant diffusion coefficient for OV} &  \text{\cite{camara2013}}   \\
\eta_i & 0.0285     & \text{Infected migrating cancer cell haptotaxis coefficient} &\text{\cite{Alzahrani2019}}\\
\eta_v & 0.0285    & \text{OV haptotaxis coefficient} & \text{\cite{Alsisi2020} } \\
\mu_1  &      0.5    & \text{Proliferation rate  for uninfected migrating cancer cells} & \text{\cite{Chaplain2006}}  \\
S_{cc} &  0.1& \text{Maximum rate of cell-cell adhesion strength} &     \text{\cite{Gerisch2008}}       \\
S_{ii} &  0.1& \text{Maximum rate of cell-cell adhesion strength} &     \text{\cite{aalsisi21}}       \\
S_{ci} &  0& \text{Maximum rate of cell-cell cross adhesion strength} &  \text{Estimated}   \\
S_{ic} &  0& \text{Maximum rate of cell-cell cross adhesion strength} &   \text{Estimated}         \\
S_{ce}      &         0.5 &  \text{Rate of Cell-ECM adhesion strength}  & \text{\cite{Painter2010}}         \\
S_{ie}      &        0.5 &  \text{Rate of Cell-ECM adhesion strength}  & \text{\cite{aalsisi21}}\\
S_{cF}      &         0.2 &  \text{Rate of Cell-fibre-ECM adhesion strength}  &  \text{\cite{Shuttleworth2019}}\\
S_{iF}      &         0.2 &  \text{Rate of Cell-fibre-ECM adhesion strength}  &  \text{Estimated}\\
\alpha_{c}         &    0.15 & \text{ECM degradation rate by uninfected cancer cells }        & \text{\cite{Alzahrani2019}}                      \\
\alpha_{i_m}       &    0.075 & \text{ECM degradation rate by infected cancer cells}      & \text{\cite{Alzahrani2019}}                         \\
\alpha_{cF}         &   0.75 & \text{Macroscopic fibre degradation rate by $c$ cells}        & \text{Estimated}\\
\alpha_{iF}         &   0.75 & \text{Macroscopic fibre degradation rate by $i$ cells}        & \text{Estimated}\\
\mu_2              &          0 & \text{Remodelling term coefficient }      &   \text{\cite{Shuttleworth2019}}  \\
\varrho       &            0.079 & \text{Infection rate of $c$ cells by OV}               & \text{\cite{Alzahrani2019}}  \\
\delta_{i}     &            0.05 &  \text{Death rate of infected cancer cells}     &   \text{\cite{camara2013}}\\
b              &      20    &   \text{Replicating rate of OVs in infected cancer cells $c$}   &  \text{\cite{camara2013}} \\
\delta_v     &   0.025& \text{Death rate of OV}   &  \text{\cite{camara2013}}\\
\nu_e         &    1  & \text{The fraction of physical space occupied by the ECM} &   \text{\cite{Shuttleworth2019}}\\
\nu_c         &    1   & \text{The fraction of physical space occupied by cancer cells} &  \text{\cite{Shuttleworth2019}}\\
\gamma_{c} & 1& \text{MDEs secretion rate by uninfected cancer cell} & \text{\cite{Shuttleworth2020}}\\
\gamma_{i} & 1.5& \text{MDEs secretion rate by infected cancer cell} & \text{\cite{Shuttleworth2020}}\\
D_m       & 0.0025 & \text{MDE diffusion coefficient}         &   \text{\cite{Peng2017}}\\
\mathcal{R}_F & 20\%:80\% &\text{The ratio of fibres and non-fibres components of ECM} & \text{\cite{Shuttleworth2019}}\\
\hline                    
\end{array}
\]
\end{table}

\subsection{\re{Baseline dynamics for the model with local flux of infected cells}}\label{res:baseline3}
\begin{figure}
    \centering
    \includegraphics[width=0.8\textwidth]{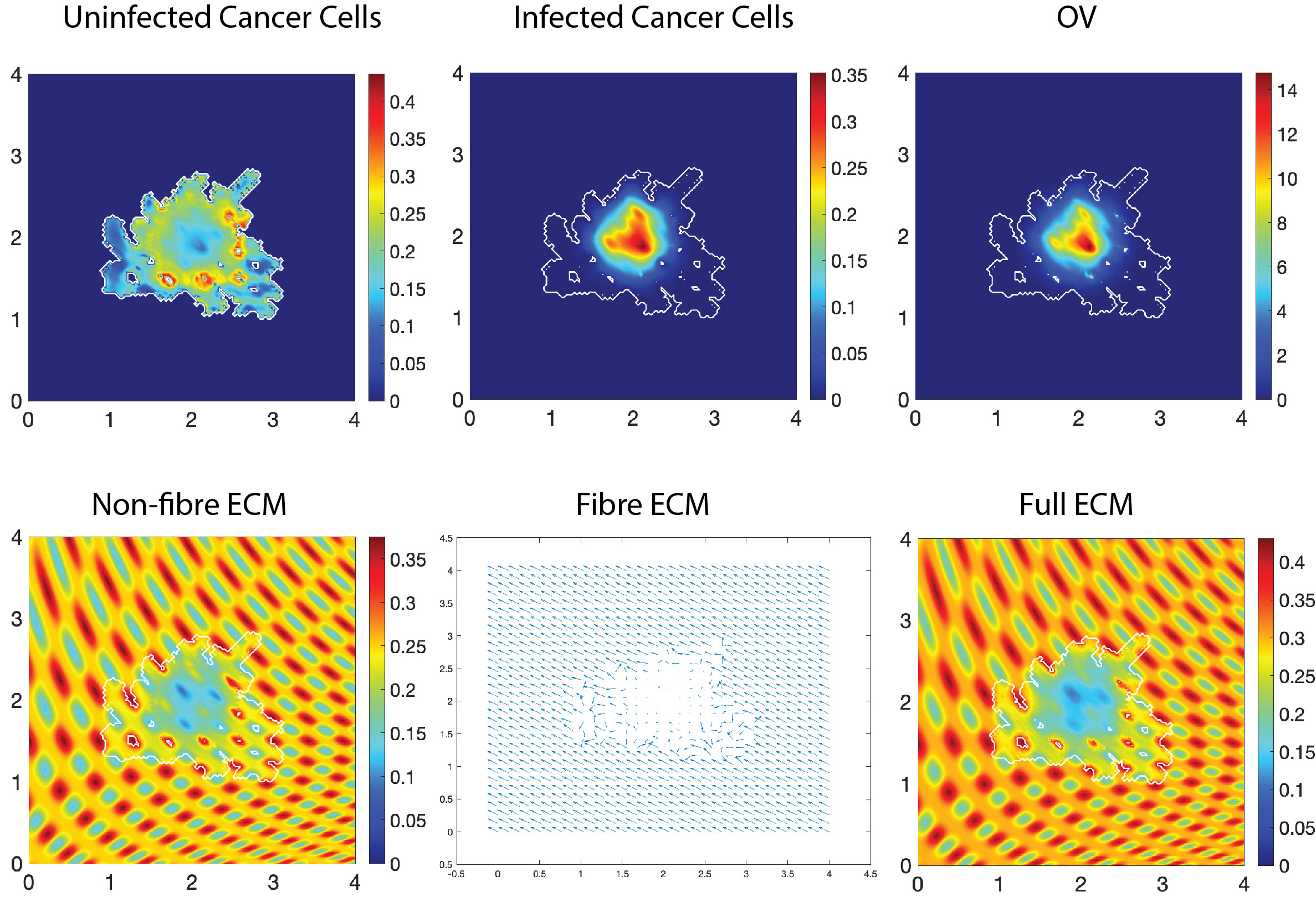}
    \caption{Simulations of system (\ref{local_sys_3}) using the parameters in Table~\ref{table_test3}. Here we show \re{the distributions of uninfected cancer cells ($c$), infected cancer cells ($i$), OV ($v$), the non-fibre ECM ($E$), fibre ECM ($F$), and the full ECM ($e$)} at micro-macro stage $75\Delta t$. }
    \label{baseline_res3}
\end{figure}
In Figure \ref{baseline_res3} \re{we present the numerical simulations obtained for the system (\ref{local_sys_3}), with} the baseline parameters in Table~\ref{table_test3}. We show the distribution of macroscopic variables at time $75 \Delta t$, \re{under the assumption that the directed movement of} infected cancer cells \re{is described by a local flux term. These simulation} results \re{also show the} distribution of non-fibre ECM and full ECM, which \re{will be omitted in the next figures} due to a lack of visibility of significant changes; only results of the oriented fibre-ECM field will be shown \re{in the next sections}.   

\subsection{The effect of \re{varying} the amount \re{of} fibres in \re{the} ECM}\label{res:fibre_amount}
\begin{figure}
    \centering
    \includegraphics[width=0.8\textwidth]{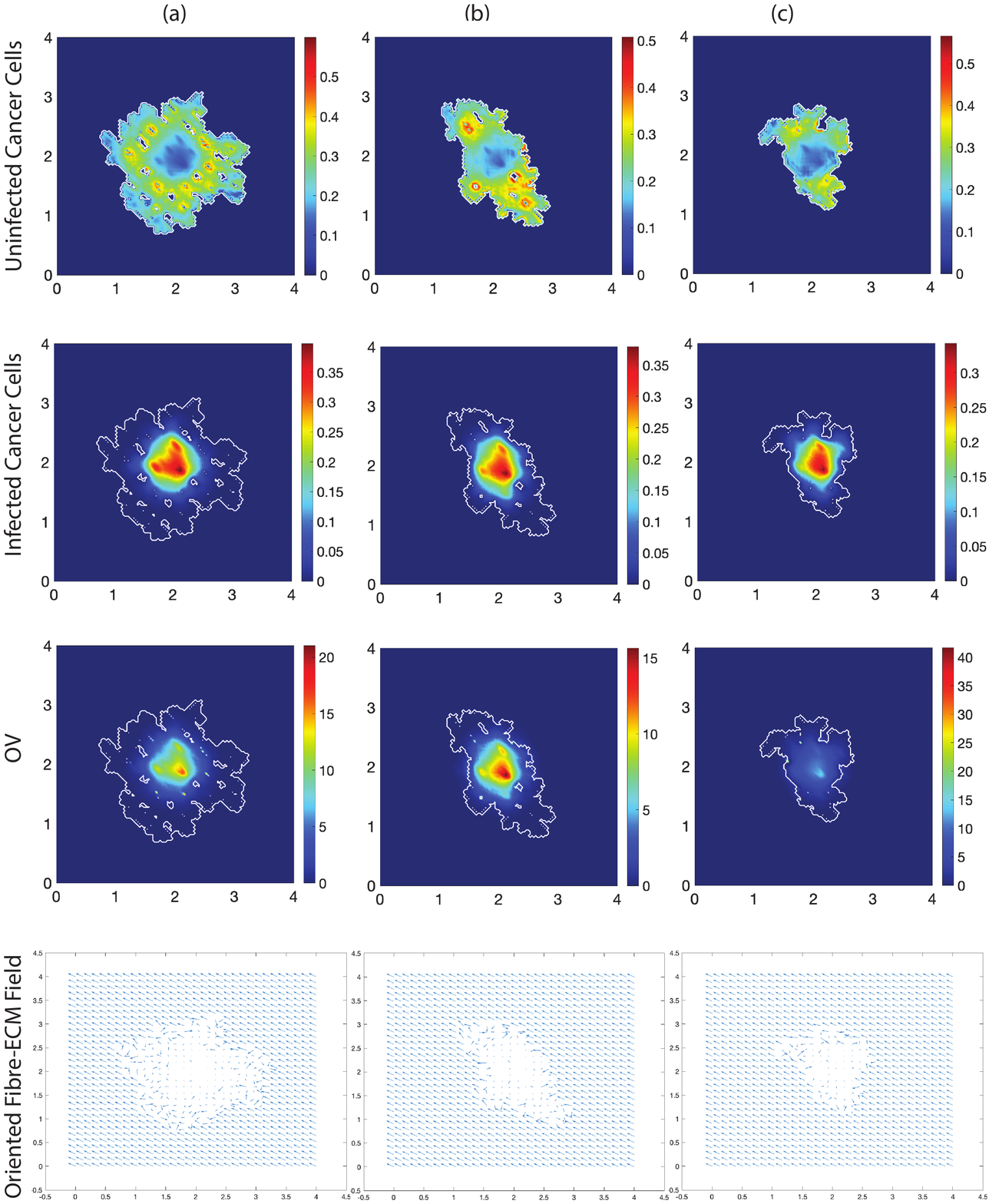}
    \caption{Simulations of system (\ref{local_sys_3}) using the parameters in Table~\ref{table_test3}. Here we show cell ($c$, $i$) and virus ($v$) distributions, \re{as well as the vector field for the oriented fibre ECM ($F$), at the micro-macro stage $75\Delta t$. We vary the ratio $\mathcal{R}_{F}$ of fibres to non-fibres components of ECM:} (a)$\mathcal{R}_F=30\%:70\%$, (b) $\mathcal{R}_F=35\%:65\%$, and (c) $\mathcal{R}_F=40\%:60\%$}
    \label{amount_fibre}
\end{figure}

\re{In} Figure \ref{amount_fibre} \re{we show} simulations of system (\ref{local_sys_3}) using the parameters in Table~\ref{table_test3}, \re{when we} vary the ratio $\mathcal{R}_F$ of fibres and non-fibres components of ECM. \re{In} Figure \ref{amount_fibre}(a) (first column) \re{we consider} $\mathcal{R}_F=30\%:70\%$, \re{in} Figure \ref{amount_fibre}(b) (second column) \re{we consider} $\mathcal{R}_F=35\%:65\%$, \re{and in} Figure \ref{amount_fibre}(c) (third column) \re{we consider} $\mathcal{R}_F=40\%:60\%$. \re{We note that} as we increase the ratio $\mathcal{R}_F$, the amount of fibres increases forcing the uninfected cancer cells to migrate \re{away} from the centre \re{of the domain, where they were initially located. This also leads to} higher OV \re{and overall lower tumour spread}. 

\subsection{The effect of increasing cell-fibre ECM adhesion strengths}\label{res:adhesion_strength3}
\begin{figure}
    \centering
    \includegraphics[width=0.8\textwidth]{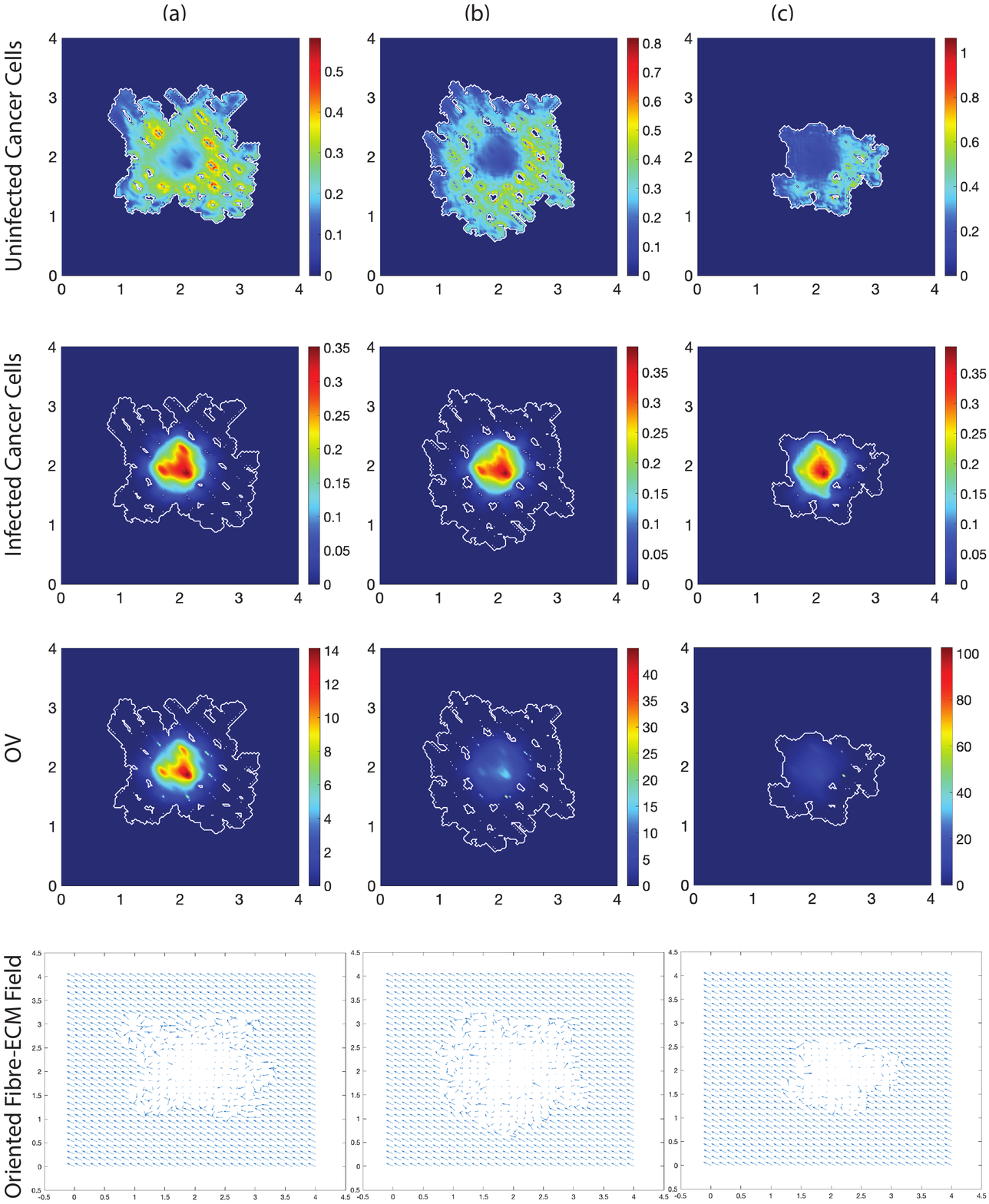}
    \caption{Simulations of system (\ref{local_sys_3}) using the parameters in Table~\ref{table_test3}. Here we show cells \re{($c$, $i$) and virus ($v$) distributions, as well as the vector field for the oriented fibre ECM ($F$) at micro-macro stage $75\Delta t$. We vary the ratio $\mathcal{R}_{F}$ of fibres to non-fibres components of ECM:} (a) $\mathcal{R}_F = 20\%:80\%$ with $S_{cF} = 0.5$, (b) $\mathcal{R}_F=30\%:70\%$ with $S_{cF} = 0.5$, and (c) $\mathcal{R}_F=40\%:60\%$ with $S_{cF} = 0.5$.}
    \label{high_fibre_ad}
\end{figure}

\re{In} Figure \ref{high_fibre_ad} \re{we show} simulations of system (\ref{local_sys_3}) \re{for} the parameters in Table~\ref{table_test3}, \re{but with} different $\mathcal{R}_F$ and \re{different} cell-fibre ECM adhesion strengths. \re{In} Figure \ref{high_fibre_ad}(a) (first column) \re{we consider} $\mathcal{R}_F=20\%:80\%$, \re{in} Figure \ref{high_fibre_ad}(b) (second column) \re{we consider} $\mathcal{R}_F=30\%:70\%$, and \re{in} Figure \ref{high_fibre_ad}(c) (third column) \re{we consider} $\mathcal{R}_F=40\%:60\%$. \re{Moreover, in all these sub-panels we take $S_{cF} = 0.5$ (compared to $S_{cF} = 0.2$ in Figure \ref{amount_fibre}). We note that} increasing cell-fibre ECM adhesion strength  \re{leads to spatial pockets of very high cancer density (i.e., $\max \;c(x, 75\Delta t)=1$ in Figure \ref{high_fibre_ad}(c), versus  $\max \;c(x, 75\Delta t)=0.55$ in Figure \ref{amount_fibre}), and even a better spatial cancer spread.}

\subsection{\re{Dynamics of the model with nonlocal flux of infected cells}}\label{res:2nonlocal}
\begin{figure}[!hp]
    \centering
    \includegraphics[width=0.8\textwidth]{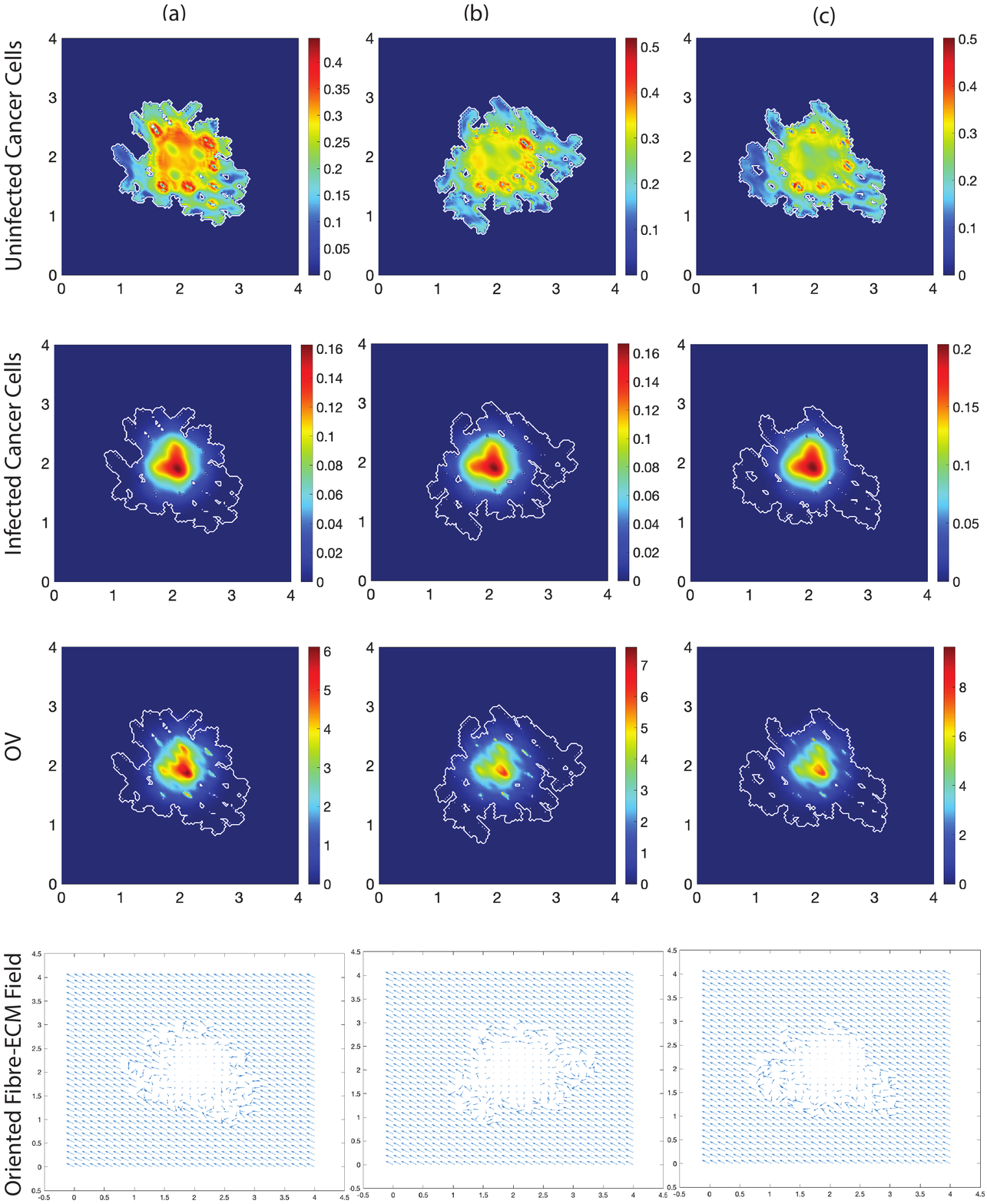}
    \caption{Simulations of system (\ref{nonlocal_sys_3}) using the parameters in Table~\ref{table_test3}. Here we show \re{the distributions of uninfected cancer cells ($c$), infected cancer cells ($i$), OV ($v$), and the vector field of the oriented fibre ECM ($F$) at micro-macro stage $75\Delta t$. (a) Baseline parameters; (b) $S_{cF} = S_{iF} = 0.3$ (while keeping all other parameters at their baseline values), (c) $S_{ie} = 0.001$ (while keeping all other parameters at their baseline values). } }
    \label{2nonlocal_fibre20}
\end{figure}

\re{In this section we investigate numerically not only} the impact of non-local \re{advection fluxes for the} infected cancer cells, \re{but also the effect of varying the} cell-cell and cell-matrix adhesion strengths. \re{In} Figure \ref{2nonlocal_fibre20}(a) \re{we show the baseline dynamics of the system (\ref{nonlocal_sys_3}) (i.e., dynamics obtained with the baseline parameters in Table~\ref{table_test3}). In} Figure \ref{2nonlocal_fibre20}(b) \re{we keep most of the parameter the same, with the exception of $S_{cF} = S_{iF} = 0.3$. In} Figure \ref{2nonlocal_fibre20}(c) \re{we keep again most of the parameter the same, with the exception of} $S_{ie} =0.001$. \re{In this case we note that varying the strengths of nonlocal cell-fibre interactions for the uninfected or infected cells does not have a significant impact on tumour structure.}

\subsection{Increasing the amount of fibres for the model with nonlocal flux of infected cells} \label{res:amount_fibre_2nonloca}

\begin{figure}[!hp]
    \centering
    \hspace{3.0cm}\includegraphics[width=0.8\textwidth]{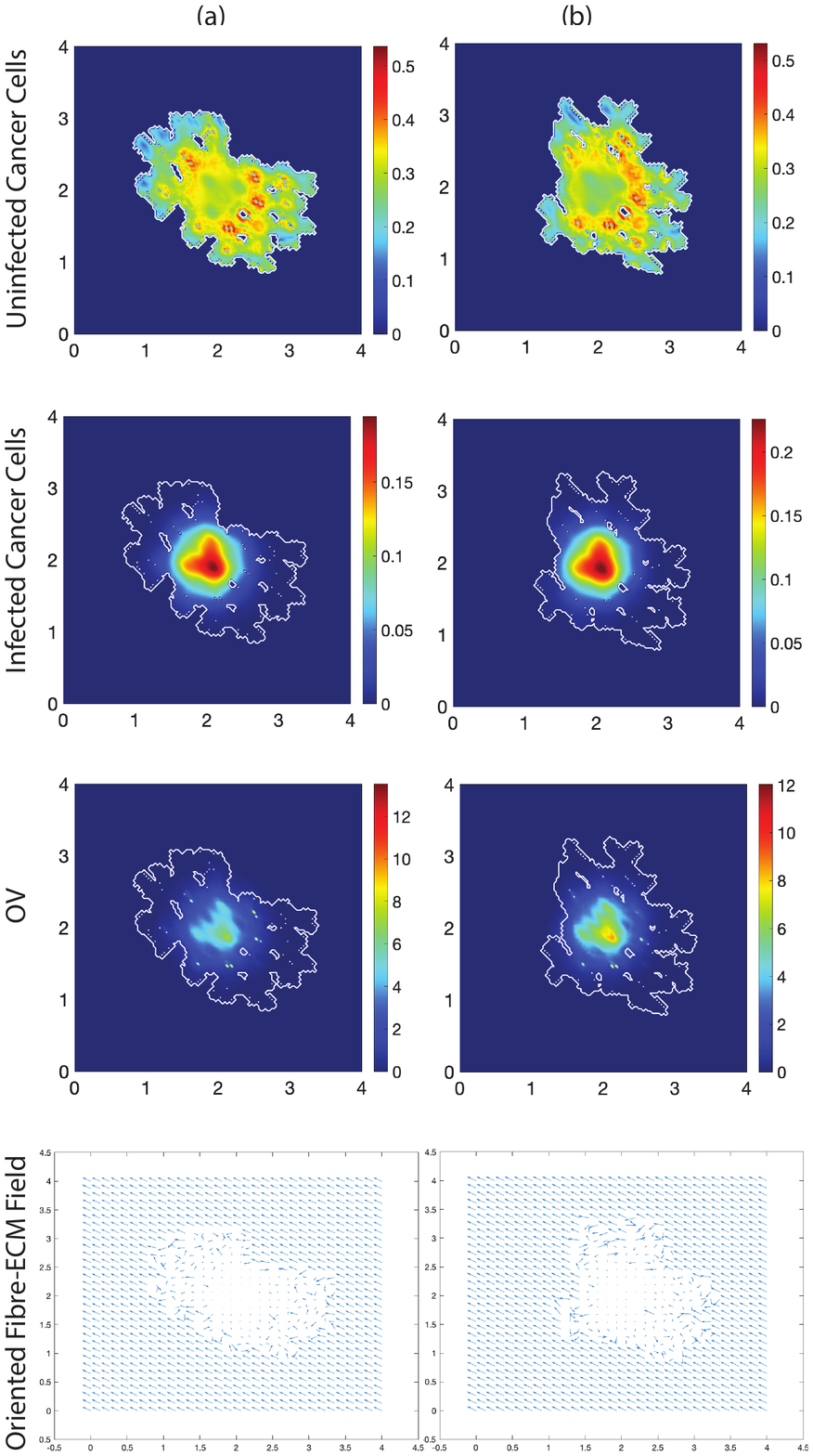}
    \caption{Simulations of system (\ref{nonlocal_sys_3}) using the parameters in Table~\ref{table_test3}. 
    \re{The sub-panels show the spatial distributions of uninfected cancer cells ($c$), infected cancer cells ($i$), OVs ($v$) and the vector field of the oriented ECM fibres ($F$) at micro-macro stage $75\Delta t$, when the ratio of fibres to non-fibres ECM components is $\mathcal{R}_F=30\%:70\%$. (a) $S_{cc}=0.1$, $S_{ce}=0.5$, (b) $S_{cc} =  0.05$, $S_{ce} =  0.001$.} }
    \label{2nonlocal_fibre30}
\end{figure}

\re{In Figure \ref{2nonlocal_fibre30} we investigate numerically the effect of increasing the amount of fibres in the ECM from $\mathcal{R}_F=20\%:80\%$ (the baseline case shown in Figure~\ref{2nonlocal_fibre20}(a)) to $\mathcal{R}_F=30\%:70\%$ (here) and varying cell-cell and cell-matrix adhesion strengths. More precisely, in sub-panels (a) we have $S_{cc}=0.1$, $S_{ce}=0.5$, while in sub-panels (b) we have $S_{cc} =  0.05$, $S_{ce} =  0.001$. The simulation results do not show significant differences between the two cases.}  

\subsection{Cross Adhesion Strength}\label{res:cross_adhesion}
\begin{figure}[!hp]
    \centering
    \hspace{3.0cm} \includegraphics[width=0.8\textwidth]{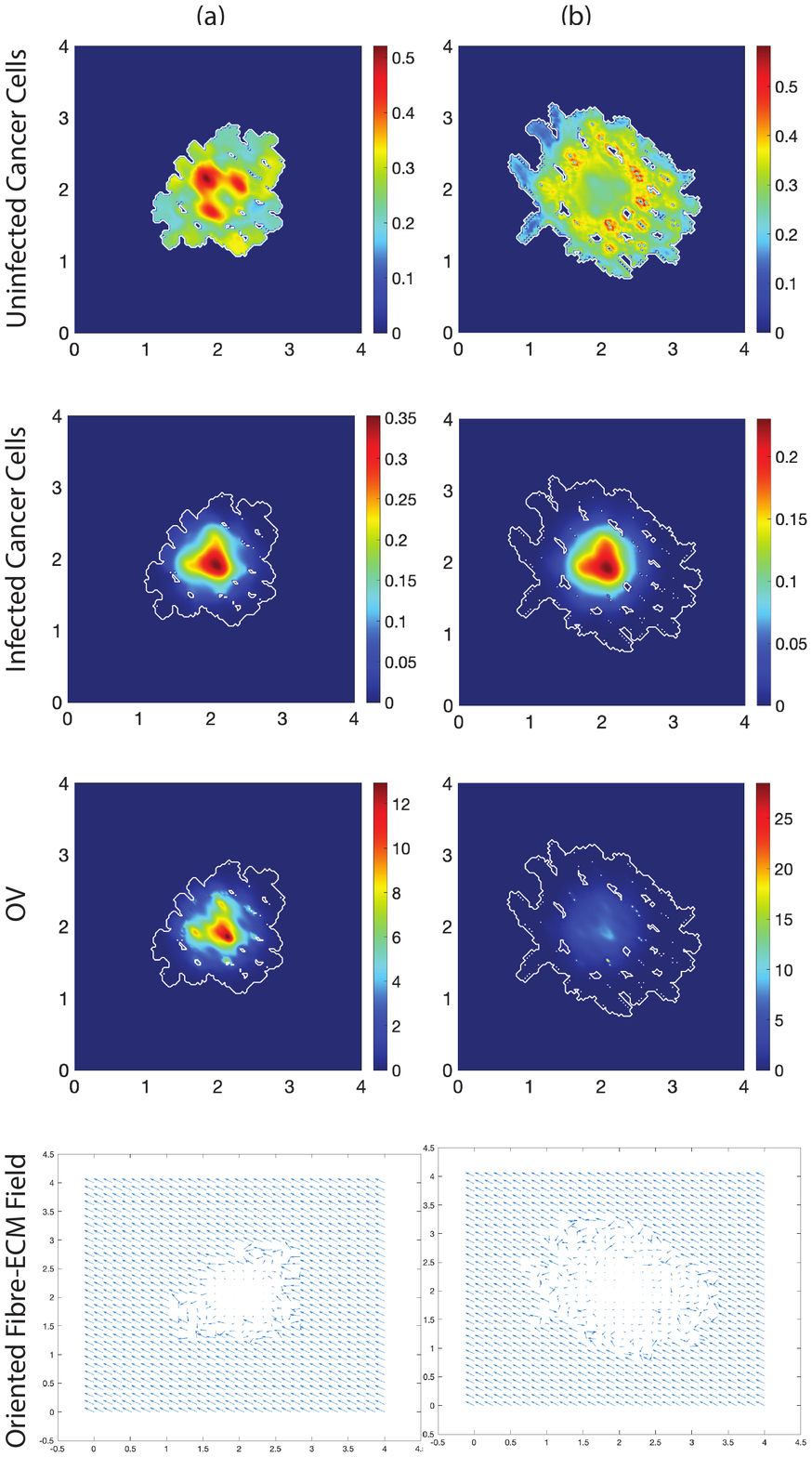}
    \caption{Simulations of system (\ref{nonlocal_sys_3}) using the parameters in Table~\ref{table_test3}. Here we show \re{the distributions of uninfected cancer cells ($c$), infected cancer cells ($i$), OVs ($v$) and the vector field of the oriented ECM fibres ($F$)} at micro-macro stage $75\Delta t$. The cell-cell and cell-matrix adhesion strengths are: (a) $S_{cc} = S_{ci} = 0.05$, $S_{ic} = S_{ii} = 0.1$ and $S_{ce} = 0.001$, (b) $S_{cc} = S_{ci} = 0.1$, $S_{ii} = S_{ic} = 0.05$ and $S_{ie} = 0.001$.}
    \label{2nonlocal_cross_ad}
\end{figure}

\re{Finally, in Figure \ref{2nonlocal_cross_ad}, we investigate numerically the effect of varying the cross adhesion strengths. In sub-panels (a) we assume that the cell-cell adhesion strengths for uninfected cancer cells ($S_{cc}$ and $S_{ci}$) are lower than the cell-cell adhesion strengths for infected cancer cells ($S_{ic}$ and $S_{ii}$). In sub-panels (b) we make the reversed assumption: the cell-cell adhesion strengths for the uninfected cells are higher than for the infected cells. We see that in this second case the tumour spreads faster through the domain (spread which is helped also by a stronger $S_{ce}=0.5$, compared to case (a) where $S_{ce}=0.001$). }  

\section{Conclusion}\label{Sect:Summary3}
In this study we \re{extended a nonlocal multi-scale moving boundary model proposed in~\cite{Alsisi2020} for oncolytic virotherapies, by considering cancer cell interactions with a heterogeneous ECM formed of fibrous and non-fibrous components. With the help of this model, we investigated numerically the impact of assumptions of local vs. non-local interactions between the infected cancer cells and uninfected cells and/or ECM. (It is known that the ECM and its components constitute a physical barrier in the spread of OVs~\cite{Wojton2010_TumEnvir-OV}, but it is not clear what are the interactions between the infected cancer cells and the environment, i.e., other cells and ECM.) }

\re{The numerical simulations showed that the ratio $\mathcal{R}_{T}$ of fibre vs. non-fibre components of the ECM, combined with the strength of cell-fibre ECM adhesion plays an important role in the extent of spatial spread of tumour cells (see Figures~\ref{amount_fibre}-\ref{high_fibre_ad}). Very large $\mathcal{R}_{T}$ ratios also seem to cause an accumulation of OVs at particular positions in space (see Figures~\ref{amount_fibre}(a)-\ref{high_fibre_ad}(c)). Moreover, the spread of the virus inside tumour seems to depend also on the strength of cell-cell interactions, with larger $S_{cc}$ causing higher viral accumulation at specific positions in space (Figure~\ref{2nonlocal_cross_ad}(b)). The nonlocal interactions between the infected cells and the environment (i.e., other cells and ECM components) play an important role in tumour and viral spread only when the magnitude of these interactions is very high (see Figure~\ref{2nonlocal_cross_ad}(a)). }

\re{This} work can be further generalised to investigate various aspects of the interactions between oncolytic virus particles and the tumour microenvironment: from the importance of \re{combining OV therapies with other classical therapies such as chemotherapies,} to the investigation of \emph{go or grow} hypothesis within a two-phase heterogeneous ECM.  

\section*{Acknowledgments}
The first author would like to acknowledge the financial support received from the Saudi Arabian Cultural Bureau in the UK on behalf of Taibah University, Medina, Saudi Arabia.

\section*{Conflict of interest}
All authors declare no conflicts of interest.

%
\bibliographystyle{AIMS}
\bibliography{library} 

\section*{Supplementary (if necessary)}

\end{document}